\documentclass[a4paper,12pt]{article}
\usepackage[top=2.5cm,bottom=2.5cm,left=2.5cm,right=2.5cm]{geometry}
\usepackage{amsmath,amssymb,amsthm,graphics,latexsym,amsfonts}
\usepackage{makecell}
\usepackage{caption}
\usepackage{booktabs}
\usepackage{graphicx}
\usepackage{bm,bbm}
\usepackage{mathrsfs}
\usepackage{tabularx}
\usepackage{cite}
\usepackage{indentfirst}
\usepackage{hhline}
\usepackage{longtable}

\newtheorem{thm}{Theorem}[section]
\newtheorem{lem}[thm]{Lemma}

\makeatletter
\newcommand{\Rmnum}[1]{\expandafter\@slowromancap\romannumeral #1@} 
\makeatother
\allowdisplaybreaks[3] 

\begin{document}

\title{{Matching forcing polynomial of generalized Petersen graph $GP(n,2)$}\thanks{Supported by the Scientific Research Project in Lanzhou University of Finance and Economics (Lzufe2019B-009), Gansu Provincial Department of Education: Youth Doctoral Fund Project (2021QB-090), Northwest China Financial Research Center Project in Lanzhou University of Finance and Economics (JYYZ201905), and Natural Science Foundation of Gansu Province (17JR5RA179).}}
\author{Shuang Zhao$^{1,2}$\\
{\small $^{1}$School of Information Engineering, Lanzhou University of Finance and Economics, }\\
{\small Lanzhou, Gansu 730000, P. R. China;}\\
{\small $^{2}$Key Laboratory of E-business Technology and Application in Gansu Province,}\\
{\small Lanzhou, Gansu 730000, P. R. China;}\\
{\small E-mail: zhaosh2018@126.com}}
\date{}
\maketitle

\begin{abstract}
    Harary et al. and Klein and Randi\'{c} proposed the forcing number of a perfect matching in mathematics and chemistry, respectively. In detail, the forcing number of a perfect matching $M$ of a graph $G$ is the smallest cardinality of subsets of $M$ that are contained in no other perfect matchings of $G$. The author and cooperators defined the forcing polynomial of $G$ as the count polynomial for perfect matchings with the same forcing number of $G$, from which the average forcing number, forcing spectrum, and the maximum and minimum forcing numbers of $G$ can be obtained. Up to now, a few papers have been considered on matching forcing problem of non-plane non-bipartite graphs. In this paper, we investigate the forcing polynomials of generalized Petersen graphs $GP(n,2)$ for $n=5,6,\ldots,15$, which is a typical class of non-plane non-bipartite graph.
\end{abstract}

\section{Introduction}
    The forcing number of a perfect matching of a graph has been introduced by Harary et al. \cite{original} in 1991. This concept can be found in earlier literatures by Klein and Randi\'{c} \cite{early,RK} under the name `innate degree of freedom', which plays an important role in the resonance theory of theoretic chemistry.

\begin{figure}[htbp]
        \centering
        \includegraphics[height=0.6in]{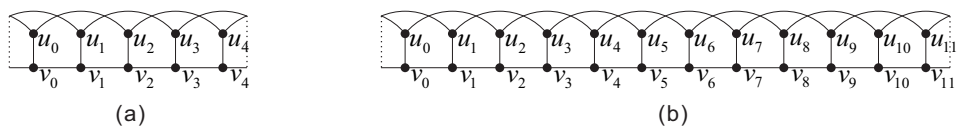}\\
        \caption{(a) The Petersen graph $GP(5,2)$ and (b) $GP(12,2)$.}
        \label{eg}
    \end{figure}

    The Petersen graph plays a special role in graph theory. In fact, the book \cite{petersen} is devoted entirely to properties and topics of the Petersen graph. Watkins \cite{def} developed the Petersen graph to try to mimic a basic structural property of it, which is known as generalized Petersen graph. In detail, a \emph{generalized Petersen graph} $GP(n,k)$ ($n\geqslant 5$, $1\leqslant k\leqslant n-1$) is a simple graph with vertex set $\{ { u }_{ i },{ v }_{ i }:0\le i\le n-1 \}$, and edge set $\{ { u }_{ i }{ u }_{ i+k },{ u }_{ i }{ v }_{ i },{ v }_{ i }{ v }_{ i+1 }:0\le i\le n-1\}$, where the subscript arithmetic is done modulo $n$ using the residues $0,1,\ldots,n-1$. Hence $GP(5,2)$ is the Petersen graph. In this paper, we concentrate on $k=2$, and place $GP(n,2)$ in a strip with the left side and right side identified as illustrated in Fig. \ref{eg}.

    Let $G$ be a graph with vertex set $V(G)$ and edge set $E(G)$. A \emph{perfect matching} $M$ of $G$ is a set of disjoint edges that covers all vertices of $G$. A cycle of $G$ is called \emph{$M$-alternating} if its edges appear alternately in $M$ and $E(G)\setminus M$. A \emph{forcing set} $S$ of $M$ is a subset of $M$ such that $S$ is contained in no other perfect matchings of $G$. The \emph{forcing number of $M$}, denoted by $f(G,M)$, is the smallest cardinality over all forcing sets of $M$. The \emph{maximum} (resp. \emph{minimum}) \emph{forcing number} of a graph $G$ is the maximum (resp. minimum) value of $f(G,M)$ over all perfect matchings $M$ of $G$, denoted by $F(G)$ (resp. $f(G)$). The \emph{forcing spectrum} of $G$ is the set of forcing numbers of all perfect matchings in $G$. In general, to compute the forcing number of a perfect matching of a graph is an NP-complete problem \cite{1}. However, some special graphs have fast algorithms to compute, such as plane bipartite graphs. The following result demonstrates the connection between forcing sets and $M$-alternating cycles in $G$.

    \begin{thm}{\em\cite{1,15}}
        \label{lem5}
        Let $G$ be a graph with a perfect matching $M.$ Then a subset $S\subseteq M$ is a forcing set of $M$ if and only if each $M$-alternating cycle of $G$ contains at least one edge of $S.$
    \end{thm}

    From the theorem we can see that the forcing number $f(G,M)$ is bounded below by the maximum number of disjoint $M$-alternating cycles, denoted by $C(G,M)$. A subgraph $H$ of a graph $G$ is said to be \emph{nice} if $G-V(H)$ has a perfect matching. By the minimax theorem on transversal, Guenin and Thomas \cite{robin} obtained the following result.

    \begin{thm}{\em\cite{robin}}
        \label{lemmima}
        Let $G$ be a bipartite graph which contains no even subdivision of $K_{3,3}$ or the Heawood graph as a nice subgraph$.$ Then for each perfect matching $M$ of $G,$ $f(G,M)=C(G,M).$
    \end{thm}

    The above theorem is a sufficient but not necessary condition on $f(G,M)=C(G,M)$ for each perfect matching $M$ of $G$. In particular, the conclusion does not hold for generalized Petersen graphs, which can be seen from the following section. Up to now, a lot of research papers talked about bipartite graphs rather than non-bipartite graphs, since in some sense it is harder to calculate non-bipartite ones.

    For a hexagonal system with a perfect matching, Xu et al. \cite{Xu} showed that the maximum forcing number is equal to the Clar number, which can measure the stability of benzenoid hydrocarbons. Also, some similar results can be found in polyomino graphs \cite{78} and (4,6)-fullerene graphs \cite{50}. Furthermore, the maximum forcing numbers of some product graphs have been studied, such as rectangle grids $P_{m}\times P_{n}$ \cite{2}, cylindrical grids $P_{m}\times C_{n}$ \cite{2,jiang2}, and tori $C_{2m}\times C_{2n}$ \cite{9}. Riddle \cite{15} derived the minimum forcing number of hypercubes $Q_k$ for even $k$ by the trailing vertex method, and Diwan \cite{diwan} obtained this for odd $k$ by algebraic method. What's more, the minimum forcing numbers of tori $C_{2m}\times C_{2n}$ \cite{15} and toroidal polyhexes \cite{yedong} were solved. Sharp lower bounds for the minimum forcing numbers of (3,6)-fullerene graphs \cite{51}, (4,6)-fullerene graphs \cite{jiangxiao} and (5,6)-fullerene graphs \cite{22,fullfor} were revealed. Moreover, the forcing spectrum of square grids $P_{2n}\times P_{2n}$ \cite{2} was also obtained. For more details, we refer the reader to a survey \cite{chechen}.

    Zhang, the author, and Lin \cite{zhao} introduced the \emph{forcing polynomial} of a graph $G$ as
    \begin{align}
        \label{equ1}
        F(G,x)=\sum_{M\in\mathcal{M}(G)}{{x}^{f(G,M)}}=\sum_{i=f(G)}^{F(G)}{\omega(G,i){x}^{i}},
    \end{align}
    where $\mathcal{M}(G)$ denotes the set of all perfect matchings of $G$, and $\omega(G,i)$ denotes the number of perfect matchings of $G$ with forcing number $i$. The following result shows that $F(G,x)$ can produce the perfect matching count, the average forcing number per perfect matching, the forcing spectrum, and the maximum and minimum forcing numbers of $G$.

    \begin{lem}{\em\cite{zhao}}
        \label{forc-perf-sum}
        The forcing polynomial of a graph $G$ has the following properties$:$\\
        $(1)$ $F(G,1)$ equals the perfect matching count of $G,$\\
        $(2)$ $\left.\frac{\mathrm{d}}{\mathrm{d}x}F(G,x)\right|_{x=1}/ F(G,1)$ equals the average forcing number per perfect matching of $G,$\\
        $(3)$ the set of degrees of $F(G,x)$ is the forcing spectrum of $G,$\\
        $(4)$ the maximum (resp. minimum) degree of $F(G,x)$ equals $F(G)$ (resp. $f(G)$).
    \end{lem}

    Afterwards, forcing polynomials of catacondensed hexagonal systems \cite{zhao}, benzenoid parallelograms \cite{zhao1}, rectangle grids \cite{forc-anti-grids}, and (5,6)-fullerene graphs C$_{60}$ \cite{buchong17}, C$_{70}$ \cite{me11} and C$_{72}$ \cite{c72} were obtained. In this paper, by definition of forcing number of a perfect matching, we show the forcing polynomials of generalized Petersen graphs $GP(n,2)$ for $n=5,6,\ldots,15$.

\section{Forcing polynomial of $GP(n,2)$}
    Up to now, there is no better method of calculating the forcing number of a perfect matching $M$ of a generalized Petersen graph $GP(n,2)$ than using the definition directly. The algorithm is to find the minimum value $k$, such that there are $k$ $M$-matched edges but containing in no other perfect matchings of $GP(n,2)$. Thus the key step is to test whether a given subset of a perfect matching is contained in other perfect matchings or not, and we call it K-step. In detail, the algorithm is given as follows: Step 1 is to initialize $k=1$. Step 2 is to pick a set of $k$ $M$-matched edges that has not been done the K-step before, and do the K-step. Step 3 is to make a decision according to different answers of K-step. If the answer is no, then stop the algorithm since we have already find the minimum forcing set of $M$. If the answer is yes and every set of $k$ $M$-matched edges has been done K-step already, then update $k$ by $k+1$ and return to Step 2. Otherwise, pick a set of $k$ $M$-matched edges that has not been done the K-step before, and return to Step 3. At the end of this algorithm, we could find the forcing number $k$ of $M$, since $k$ can not exceed $\frac{|V(GP(n,2))|}{2}$, namely the number of edges in $M$.

    As an example we calculate the forcing polynomial of the Petersen graph $GP(5,2)$. Obviously $GP(5,2)$ has six perfect matchings, which are $M_1=\{u_0u_2,u_1u_3,u_4v_4,v_0v_1,v_2v_3\}$, $M_2=\{u_1u_3,u_2u_4,u_0v_0,v_1v_2,v_3v_4\}$, $M_3=\{u_2u_4,u_3u_0,u_1v_1,v_2v_3,v_4v_0\}$, $M_4=\{u_3u_0,u_4u_1,$\\$u_2v_2,v_3v_4,v_0v_1\}$, $M_5=\{u_4u_1,u_0u_2,u_3v_3,v_4v_0,v_1v_2\}$, $M_6=\{u_0v_0,u_1v_1,u_2v_2,u_3v_3,u_4v_4\}$. There is an automorphism $g_j$ of $GP(5,2)$ such that $g_j(M_1)=M_{1+j}$ for $j=1,2,3,4$, where $g_j(u_i)=u_{i+j}$ and $g_j(v_i)=v_{i+j}$ for $i=0,1,2,3,4$. Hence $f(GP(5,2),M_1)=f(GP(5,2),M_{1+j})$, and we call $M_1$ and $M_{1+j}$ are \emph{equivalent structures} for $j=1,2,3,4$. In general, two perfect matchings $M',M''$ of $GP(n,2)$ are called \emph{equivalent structures} if there is an automorphism $g$ of $GP(n,2)$ and $j\in \{1,2,\ldots,n-1\}$ such that $g(M')=M''$, where $g(u_i)=u_{i+j}$, $g(v_i)=v_{i+j}$ for $i=0,1,\ldots,n-1$.

    First we consider the forcing number of $M_1$. Since $u_0u_2\in M_1\cap M_5$, $u_1u_3\in M_1\cap M_2$, $u_4v_4\in M_1\cap M_6$, $v_0v_1\in M_1\cap M_4$, $v_2v_3\in M_1\cap M_3$, we have $f(GP(5,2),M_1)\geqslant 2$. Since $\{u_0u_2,u_1u_3\}\subseteq M_1$, but $\{u_0u_2,u_1u_3\}\not\subseteq M_i$ for $i=2,3,4,5,6$, we have $f(GP(5,2),M_1)=2$. Next we consider the forcing number of $M_6$. Since $u_4v_4\in M_6\cap M_1$ and $u_kv_k\in M_6\cap M_{k+2}$ for $k=0,1,2,3$, and $\{u_0v_0,u_1v_1\}\subseteq M_6$ but $\{u_0v_0,u_1v_1\}\not\subseteq M_i$ for $i=1,2,3,4,5$, we have $f(GP(5,2),M_6)=2$. It follows that $F(GP(5,2),x)=6x^2$.

    In the following table, we illustrate two sets of bold lines in the first and second graphs in the column NES to represent $M_1$ and $M_6$, together with corresponding minimum forcing sets illustrated with two sets of double lines, respectively. Note that in the following tables, NES is short for non-equivalent structures of generalized Petersen graph, PMC is short for the corresponding perfect matching count of equivalent structure, FN is short for the forcing number of perfect matching, and FP-$n$ is short for the forcing polynomial of generalized Petersen graph for $n=5,6,\ldots,15$.

    \begin{longtable}{c|c|cc}
        \toprule
        NO & NES & PMC & FN  \\
        \hline
        1 & \makecell{\includegraphics[height=0.3in]{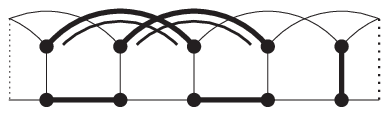}} & 5 & 2 \\
        \hline
        2 & \makecell{\includegraphics[height=0.29in]{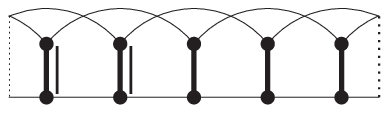}} & 1 & 2 \\
        \hline \hline
        FP-5 & \multicolumn {3}{l}{$6x^2$}\\
        \bottomrule
    \end{longtable}

    Note that $f(GP(5,2),M)=2$ for every perfect matching $M$ of $GP(5,2)$. Furthermore, it is not hard to derive that all the $M_1$-alternating cycles are $u_0u_2v_2v_3u_3u_1v_1v_0u_0$, $u_0u_2v_2v_3v_4u_4u_1u_3u_0$, $u_0u_2u_4v_4v_3v_2v_1v_0u_0$, $u_0u_2u_4v_4v_0v_1u_1u_3u_0$, $u_1u_3v_3v_2v_1v_0v_4u_4u_1$. Also, all the $M_6$-alternating cycles are $u_iv_iv_{1+i}u_{1+i}u_{3+i}v_{3+i}v_{2+i}u_{2+i}u_i$ for $i=0,1,2,3,4$. Hence $C(GP(5,2),M)=1$ for every perfect matching $M$ of $GP(5,2)$. This implies that $f(GP(5,2),M)> C(GP(5,2),M)$.

    By a similar argument to the above, we can derive the forcing polynomial of other generalized Petersen graphs $GP(n,2)$ for $n=6,7,\ldots,15$ in the following tables.

    \begin{longtable}{c|c|cc}
        \toprule
        NO & NES & PMC & FN  \\
        \hline
        1 & \makecell{\includegraphics[height=0.3in]{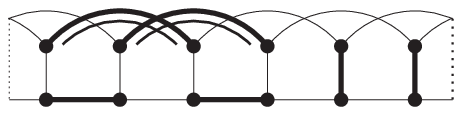}} & 6 & 2 \\
        \hline
        2 & \makecell{\includegraphics[height=0.29in]{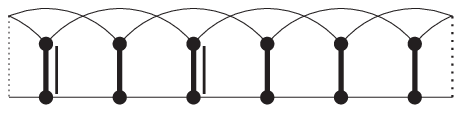}} & 1 & 2 \\
        \hline
        3 & \makecell{\includegraphics[height=0.3in]{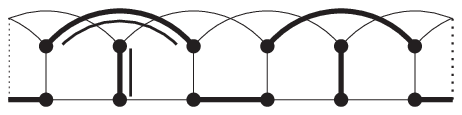}} & 3 & 2 \\
        \hline \hline
        FP-6 & \multicolumn {3}{l}{$10x^2$}\\
        \bottomrule
    \end{longtable}

    \begin{longtable}{c|c|cc}
        \toprule
        NO & NES & PMC & FN  \\
        \hline
        1 & \makecell{\includegraphics[height=0.3in]{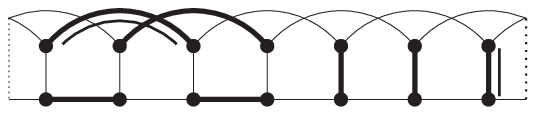}} & 7 & 2 \\
        \hline
        2 & \makecell{\includegraphics[height=0.29in]{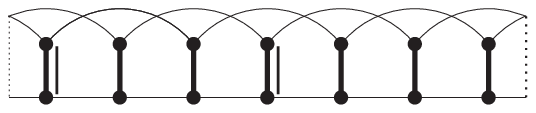}} & 1 & 2 \\
        \hline
        3 & \makecell{\includegraphics[height=0.3in]{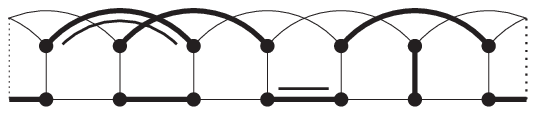}} & 7 & 2 \\
        \hline \hline
        FP-7 & \multicolumn {3}{l}{$15x^2$}\\
        \bottomrule
    \end{longtable}

    \begin{longtable}{c|c|cc}
        \toprule
        NO & NES & PMC & FN  \\
        \hline
        1 & \makecell{\includegraphics[height=0.3in]{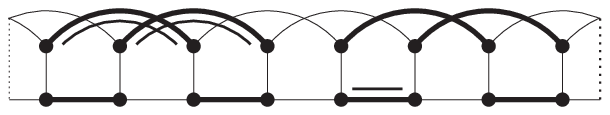}} & 4 & 3 \\
        \hline
        2 & \makecell{\includegraphics[height=0.3in]{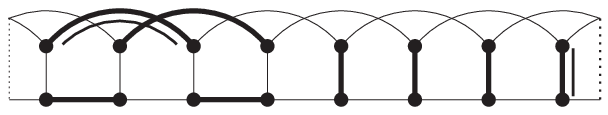}} & 8 & 2 \\
        \hline
        3 & \makecell{\includegraphics[height=0.3in]{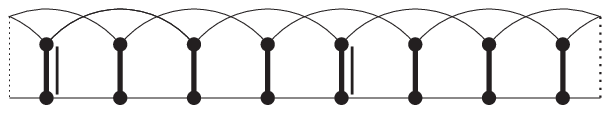}} & 1 & 2 \\
        \hline
        4 & \makecell{\includegraphics[height=0.3in]{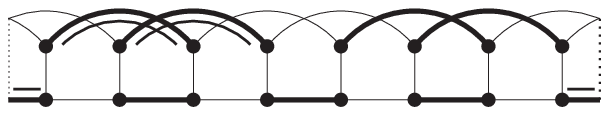}} & 4 & 3 \\
        \hline \hline
        FP-8 & \multicolumn {3}{l}{$8x^3+9x^2$}\\
        \bottomrule
    \end{longtable}

    \begin{longtable}{c|c|cc}
        \toprule
        NO & NES & PMC & FN  \\
        \hline
        1 & \makecell{\includegraphics[height=0.3in]{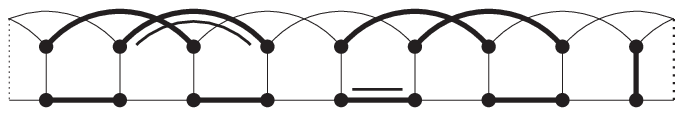}} & 9 & 2 \\
        \hline
        2 & \makecell{\includegraphics[height=0.3in]{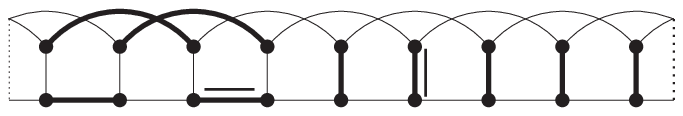}} & 9 & 2 \\
        \hline
        3 & \makecell{\includegraphics[height=0.3in]{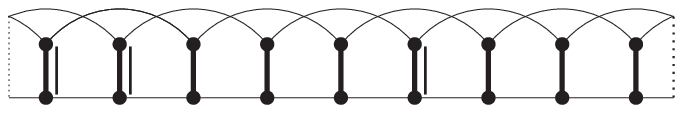}} & 1 & 3 \\
        \hline
        4 & \makecell{\includegraphics[height=0.3in]{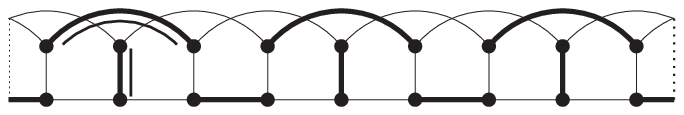}} & 3 & 2 \\
        \hline \hline
        FP-9 & \multicolumn {3}{l}{$x^3+21x^2$}\\
        \bottomrule
    \end{longtable}

    \begin{longtable}{c|c|cc}
        \toprule
        NO & NES & PMC & FN  \\
        \hline
        1 & \makecell{\includegraphics[height=0.3in]{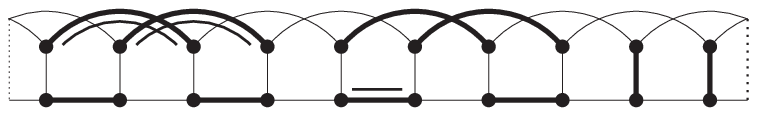}} & 10 & 3 \\
        \hline
        2 & \makecell{\includegraphics[height=0.3in]{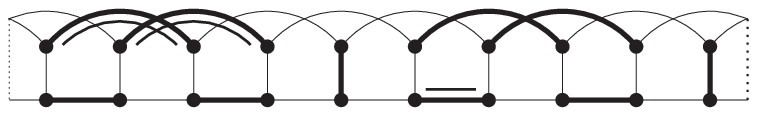}} & 5 & 3 \\
        \hline
        3 & \makecell{\includegraphics[height=0.3in]{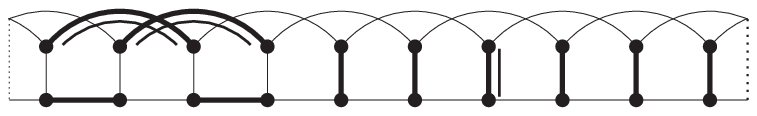}} & 10 & 3 \\
        \hline
        4 & \makecell{\includegraphics[height=0.29in]{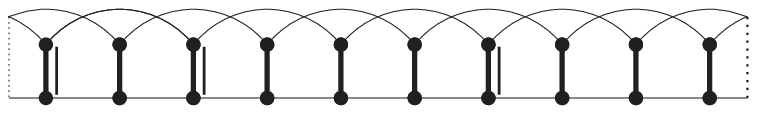}} & 1 & 3 \\
        \hline
        5 & \makecell{\includegraphics[height=0.3in]{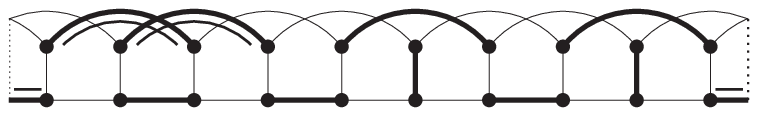}} & 10 & 3 \\
        \hline \hline
        FP-10 & \multicolumn {3}{l}{$36x^3$}\\
        \bottomrule
    \end{longtable}

    \begin{longtable}{c|c|cc}
        \toprule
        NO & NES & PMC & FN  \\
        \hline
        1 & \makecell{\includegraphics[height=0.3in]{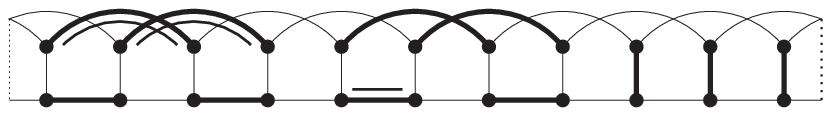}} & 11 & 3 \\
        \hline
        2 & \makecell{\includegraphics[height=0.3in]{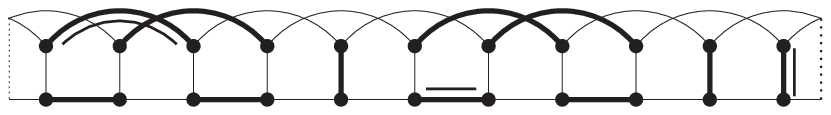}} & 11 & 3 \\
        \hline
        3 & \makecell{\includegraphics[height=0.3in]{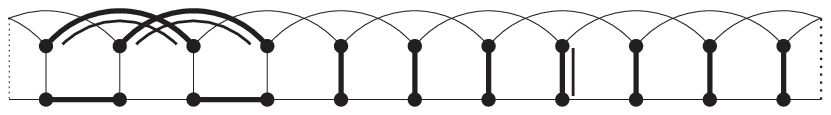}} & 11 & 3 \\
        \hline
        4 & \makecell{\includegraphics[height=0.3in]{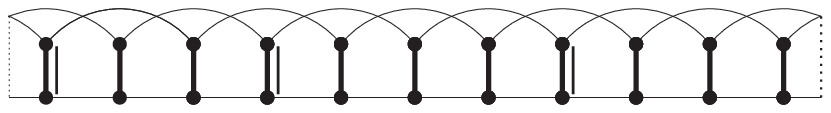}} & 1 & 3 \\
        \hline
        5 & \makecell{\includegraphics[height=0.3in]{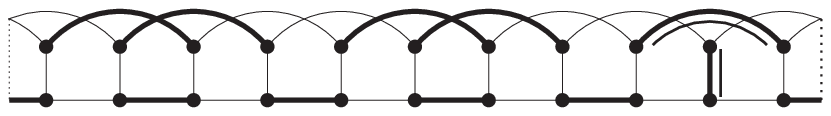}} & 11 & 2 \\
        \hline \hline
        FP-11 & \multicolumn {3}{l}{$34x^3+11x^2$}\\
        \bottomrule
    \end{longtable}

    \begin{longtable}{c|c|cc}
        \toprule
        NO & NES & PMC & FN  \\
        \hline
        1 & \makecell{\includegraphics[height=0.3in]{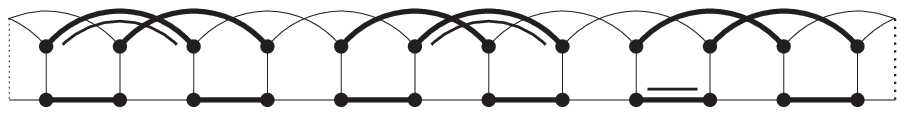}} & 4 & 3 \\
        \hline
        2 & \makecell{\includegraphics[height=0.3in]{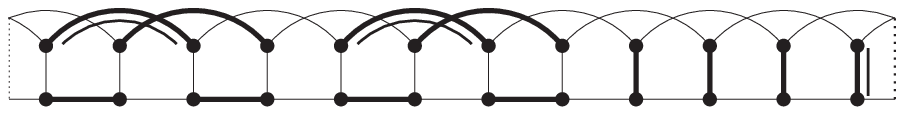}} & 12 & 3 \\
        \hline
        3 & \makecell{\includegraphics[height=0.3in]{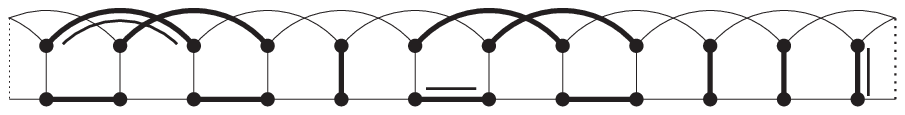}} & 12 & 3 \\
        \hline
        4 & \makecell{\includegraphics[height=0.3in]{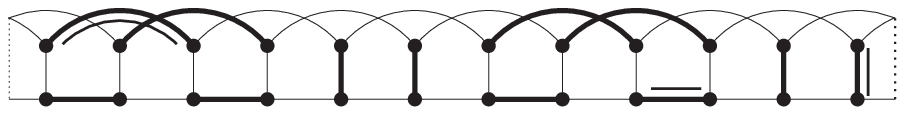}} & 6 & 3 \\
        \hline
        5 & \makecell{\includegraphics[height=0.3in]{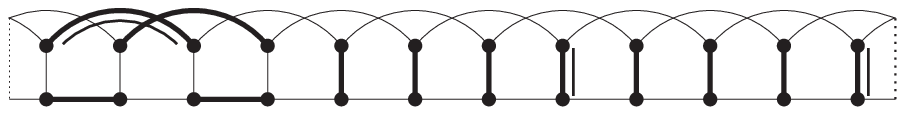}} & 12 & 3 \\
        \hline
        6 & \makecell{\includegraphics[height=0.3in]{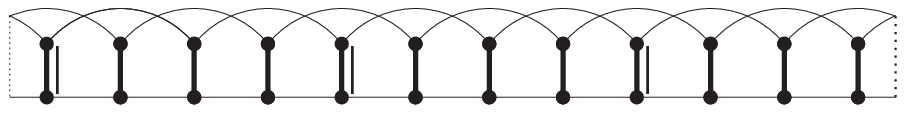}} & 1 & 3 \\
        \hline
        7 & \makecell{\includegraphics[height=0.3in]{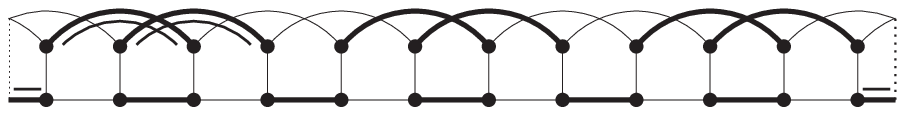}} & 4 & 3 \\
        \hline
        8 & \makecell{\includegraphics[height=0.3in]{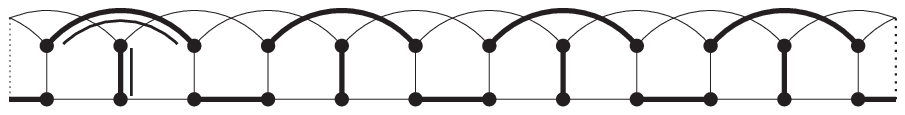}} & 3 & 2 \\
        \hline \hline
        FP-12 & \multicolumn {3}{l}{$51x^3+3x^2$}\\
        \bottomrule
    \end{longtable}

    \begin{longtable}{c|c|cc}
        \toprule
        NO & NES & PMC & FN  \\
        \hline
        1 & \makecell{\includegraphics[height=0.3in]{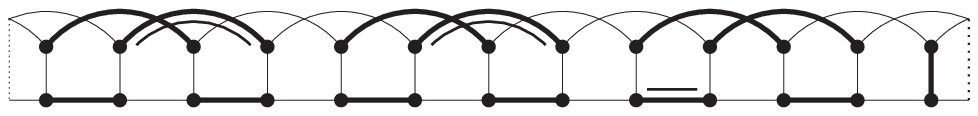}} & 13 & 3 \\
        \hline
        2 & \makecell{\includegraphics[height=0.3in]{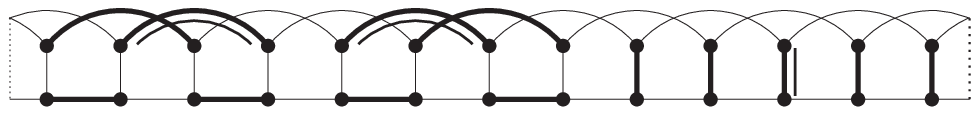}} & 13 & 3 \\
        \hline
        3 & \makecell{\includegraphics[height=0.3in]{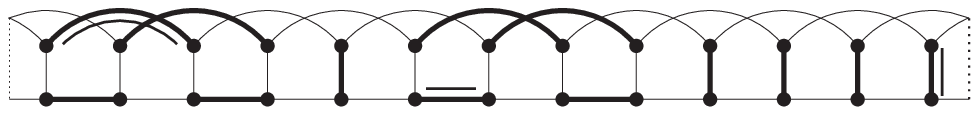}} & 13 & 3 \\
        \hline
        4 & \makecell{\includegraphics[height=0.3in]{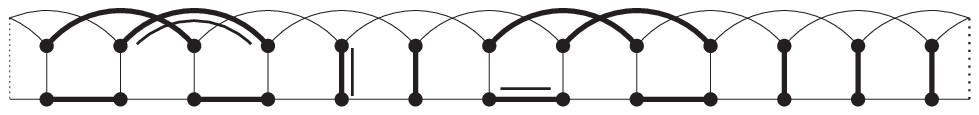}} & 13 & 3 \\
        \hline
        5 & \makecell{\includegraphics[height=0.3in]{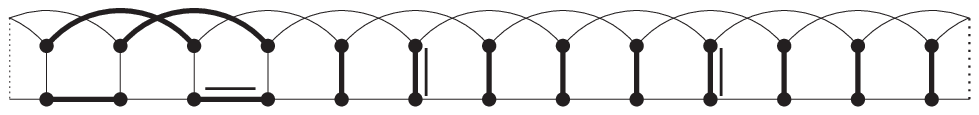}} & 13 & 3 \\
        \hline
        6 & \makecell{\includegraphics[height=0.3in]{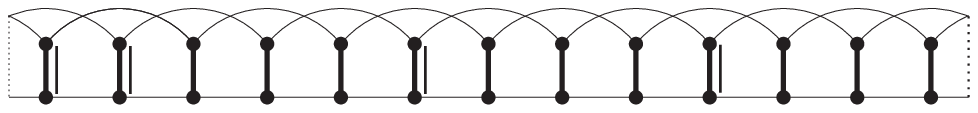}} & 1 & 4 \\
        \hline
        7 & \makecell{\includegraphics[height=0.3in]{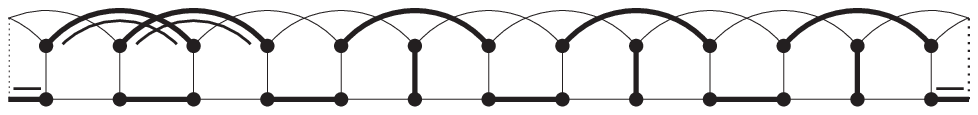}} & 13 & 3 \\
        \hline \hline
        FP-13 & \multicolumn {3}{l}{$x^4+78x^3$}\\
        \bottomrule
    \end{longtable}

    \begin{longtable}{c|c|cc}
        \toprule
        NO & NES & PMC & FN  \\
        \hline
        1 & \makecell{\includegraphics[height=0.3in]{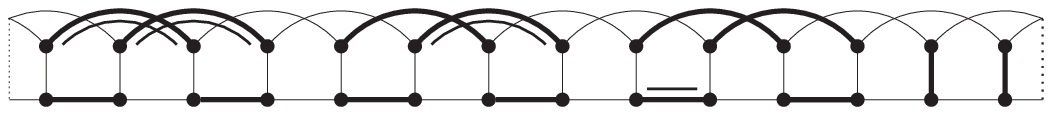}} & 14 & 4 \\
        \hline
        2 & \makecell{\includegraphics[height=0.3in]{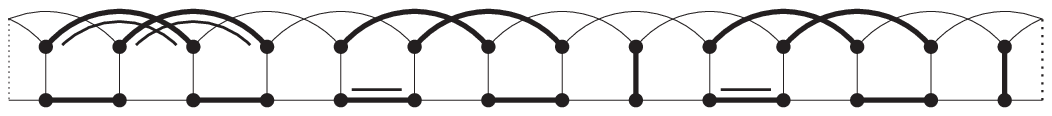}} & 14 & 4 \\
        \hline
        3 & \makecell{\includegraphics[height=0.3in]{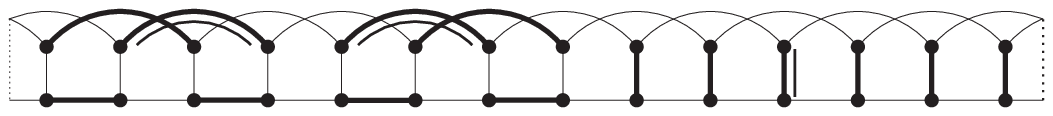}} & 14 & 3 \\
        \hline
        4 & \makecell{\includegraphics[height=0.3in]{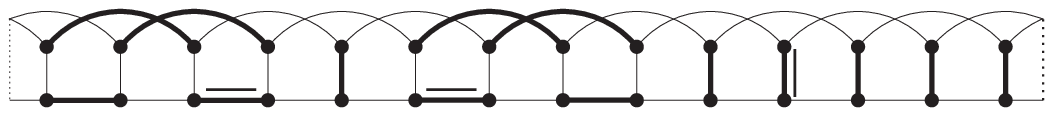}} & 14 & 3 \\
        \hline
        5 & \makecell{\includegraphics[height=0.3in]{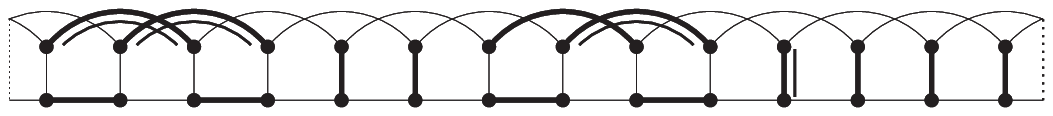}} & 14 & 4 \\
        \hline
        6 & \makecell{\includegraphics[height=0.3in]{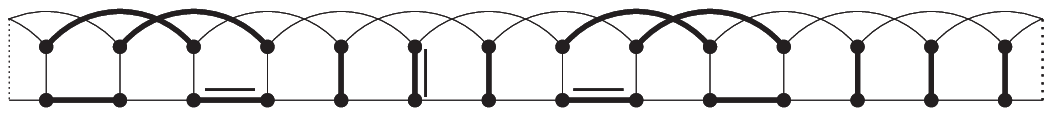}} & 7 & 3 \\
        \hline
        7 & \makecell{\includegraphics[height=0.3in]{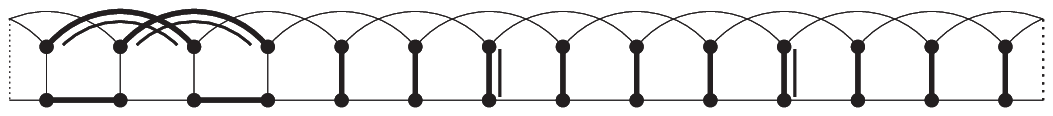}} & 14 & 4 \\
        \hline
        8 & \makecell{\includegraphics[height=0.3in]{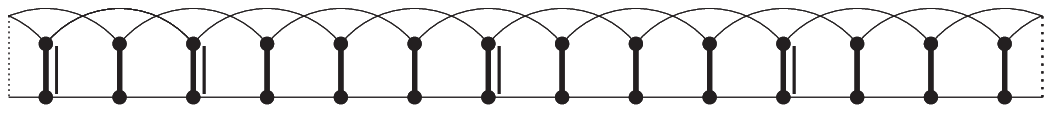}} & 1 & 4 \\
        \hline
        9 & \makecell{\includegraphics[height=0.3in]{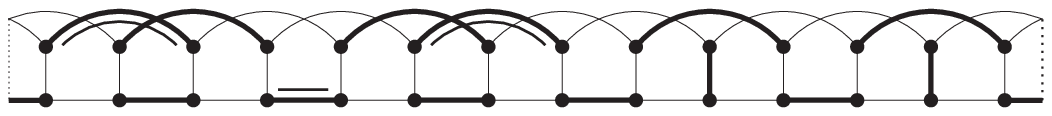}} & 14 & 3 \\
        \hline
        10 & \makecell{\includegraphics[height=0.3in]{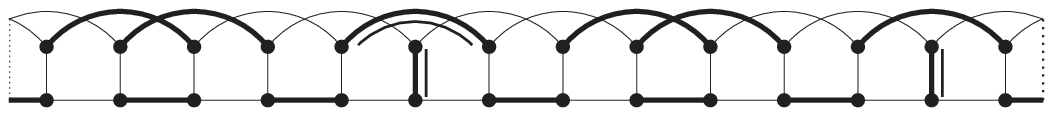}} & 7 & 3 \\
        \hline \hline
        FP-14 & \multicolumn {3}{l}{$57x^4+56x^3$}\\
        \bottomrule
    \end{longtable}

    \begin{longtable}{c|c|cc}
        \toprule
        NO & NES & PMC & FN  \\
        \hline
        1 & \makecell{\includegraphics[height=0.3in]{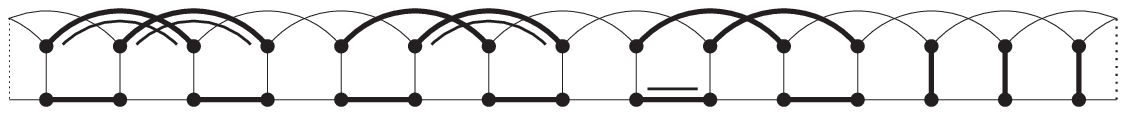}} & 15 & 4 \\
        \hline
        2 & \makecell{\includegraphics[height=0.3in]{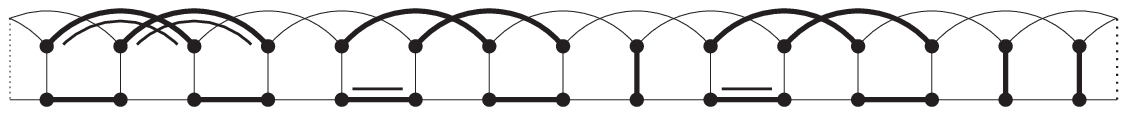}} & 15 & 4 \\
        \hline
        3 & \makecell{\includegraphics[height=0.3in]{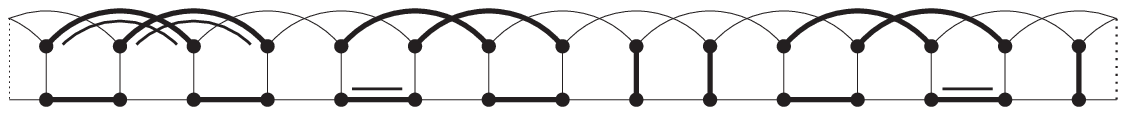}} & 15 & 4 \\
        \hline
        4 & \makecell{\includegraphics[height=0.3in]{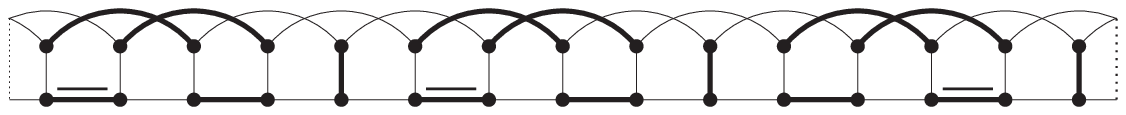}} & 5 & 3  \\
        \hline
        5 & \makecell{\includegraphics[height=0.3in]{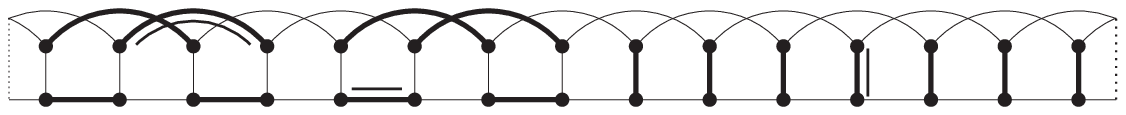}} & 15 & 3 \\
        \hline
        6 & \makecell{\includegraphics[height=0.3in]{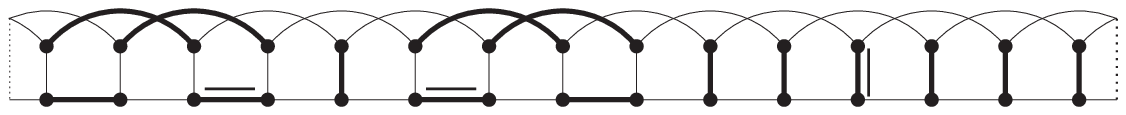}} & 15 & 3 \\
        \hline
        7 & \makecell{\includegraphics[height=0.3in]{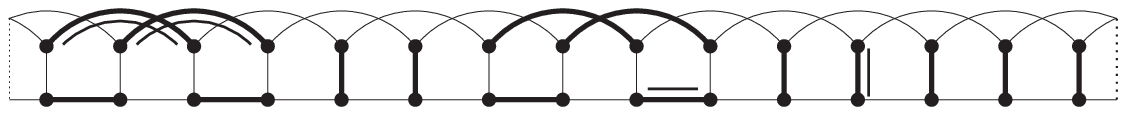}} & 15 & 4 \\
        \hline
        8 & \makecell{\includegraphics[height=0.3in]{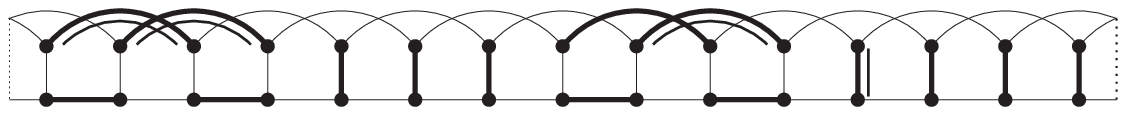}} & 15 & 4 \\
        \hline
        9 & \makecell{\includegraphics[height=0.3in]{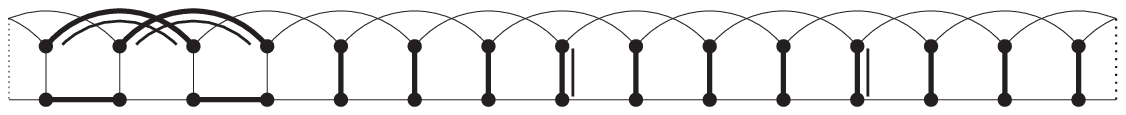}} & 15 & 4 \\
        \hline
        10 & \makecell{\includegraphics[height=0.29in]{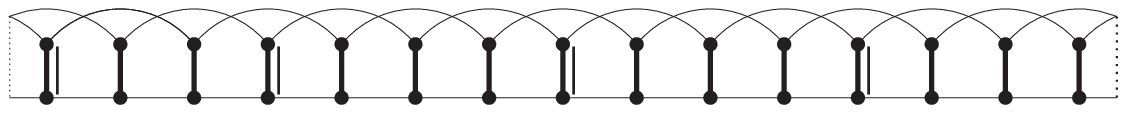}} & 1 & 4 \\
        \hline
        11 & \makecell{\includegraphics[height=0.3in]{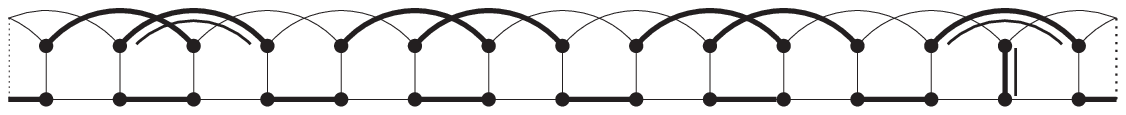}} & 15 & 3 \\
        \hline
        12 & \makecell{\includegraphics[height=0.3in]{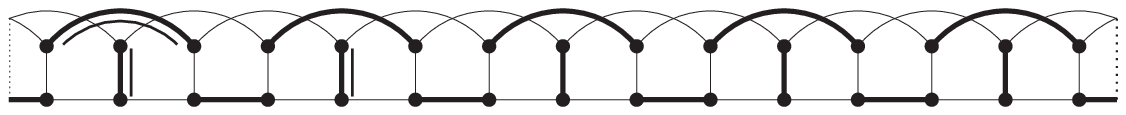}} & 3 & 3 \\
        \hline \hline
        FP-15 & \multicolumn {3}{l}{$91x^4+53x^3$}\\
        \bottomrule
    \end{longtable}

\end{document}